\documentclass[preprint,12pt]{elsarticle}

\usepackage{amssymb}
\usepackage{multirow}
\usepackage{amsthm}
\usepackage{graphics}
\usepackage{epsfig}
\usepackage{amssymb}
\usepackage{graphicx}
\usepackage{subfigure}
\usepackage{color}

\usepackage{tikz}
\usetikzlibrary{calc}
\usetikzlibrary{positioning}
\usetikzlibrary{patterns}
\usepackage{hyperref}
\usepackage{a4wide}
\usepackage{amsmath}
\usepackage{amsfonts}
\usepackage{enumerate}
\newcommand{\pf}{\noindent {\bf Proof: }}

\newtheorem*{theoremaux}{Theorem \theoremauxnum}
\gdef\theoremauxnum{1}

\newtheorem{lemma}{\bf Lemma}[section]
 
\newtheorem{theorem}{\bf Theorem}[section]
\newtheorem{proposition}[lemma]{\bf Proposition}
\newtheorem{corollary}[lemma]{\bf Corollary}
\newtheorem{definition}{\bf Definition}[section]
\newtheorem{remark}{\bf Remark}[section]

\journal{~}

\begin{document}

\begin{frontmatter}



\title{The Difference Subgroup Graph of a Finite Group}



\author{Angsuman Das\corref{cor1}}
\ead{angsuman.maths@presiuniv.ac.in}
\author{Arnab Mandal}
\ead{arnab.maths@presiuniv.ac.in}
\author{Labani Sarkar}
\ead{labanis890@gmail.com }
\address{Department of Mathematics, Presidency University, Kolkata, India}

\cortext[cor1]{Corresponding author}

\begin{abstract}
The \emph{difference subgroup graph} \(D(G)\) of a finite group \(G\) is defined as the graph whose vertices are the non-trivial proper subgroups of \(G\), with two distinct vertices \(H\) and \(K\) adjacent if and only if \(\langle H, K \rangle = G\) but \(HK \ne G\). This graph arises naturally as the difference between the join graph \(\Delta(G)\) and the comaximal subgroup graph \(\Gamma(G)\). In this paper, we initiate a systematic study of \(D(G)\) and its reduced version \(D^*(G)\), obtained by removing isolated vertices.

We establish several fundamental structural properties of these graphs, including conditions for connectivity, forbidden subgraph characterizations, and the relationship between graph parameters — such as independence number, clique number, and girth — and the solvability or nilpotency of the underlying group. 

The paper concludes with a discussion of open problems and potential directions for future research.
\end{abstract}


\end{frontmatter}

\section{Introduction}
The representation of algebraic structures as graphs is a fertile area of research, tracing its origins to Cayley graphs. The primary aim is to uncover structural properties of the algebraic object by studying the associated graph. In the context of groups, these graphs broadly fall into two categories based on their vertex sets: those whose vertices are elements of the group, and those whose vertices are subgroups.

The first class includes many important graphs such as power graphs  \cite{Cameron-Ghosh}, enhanced power graphs \cite{enhanced}, commuting graphs \cite{commuting-graph} and difference graphs \cite{difference_graph} where the vertices represent elements of the group. The second class consists of graphs like comaximal subgroup graphs \cite{akbari}, \cite{1st-paper} join graphs \cite{ahmadi-taeri}, subgroup inclusion graphs \cite{inclusion} and subgroup intersection graphs \cite{intersection-subgroup-graph}, whose vertices are the non-trivial proper subgroups of the underlying group.

While extensive literature exists on the first type of graphs (for a comprehensive survey on this topic, see \cite{Cameron-survey}), the study of the second type remains fragmented. To date, no systematic work has been done to compare and analyze different graphs in this category.

Targeting this goal, in this paper, we study the graphs of the second type. Before proceeding further, let us recall the definitions of two important graphs of this type, namely comaximal subgroup graphs and join graphs.

\begin{definition} Let $G$ be a group and $S$ be the collection of all non-trivial proper subgroups of $G$. The co-maximal subgroup graph $\Gamma(G)$ of a group $G$ is defined to be a graph with $S$ as the set of vertices and two distinct vertices $H$ and $K$ are adjacent if and only if $HK=G$.	The deleted co-maximal subgroup graph of $G$, denoted by $\Gamma^*(G)$, is defined as the graph obtained by removing the isolated vertices from $\Gamma(G)$.

The join of subgroup graph $\Delta(G)$ of a group $G$ is defined to be a graph with $S$ as the set of vertices and two distinct vertices $H$ and $K$ are adjacent if and only if $\langle H,K\rangle=G$.
\end{definition}
\begin{remark}
    It is to be noted that although $\Delta(G)$ was introduced in \cite{ahmadi-taeri}, we use a slightly modified version of it as in \cite{lucchini-join}. We take all non-trivial proper subgroups as vertices, instead of taking only those proper subgroups of $G$ which are not contained in the Frattini subgroup of $G$. This change allows us to study both the graphs, namely $\Gamma(G)$ and $\Delta(G)$, in the same platform, i.e., under this modification, $\Gamma(G)$ is a subgraph of $\Delta(G)$.
\end{remark}

As $\Gamma(G)$ is a subgraph of $\Delta(G)$, we can naturally define a graph $D(G)$ which is the difference of these two graphs, i.e., a graph which has same set of vertices as that of $\Gamma(G)$ and $\Delta(G)$, and two vertices are adjacent if they are adjacent in $\Delta(G)$ but not adjacent in $\Gamma(G)$. It is to be noted that some recent works \cite{difference_graph}, \cite{kumar-difference} have focused on the difference graph of enhanced power graphs and power graphs, two graphs from the first class as mentioned earlier.

\subsection{Our Contribution}

We begin by presenting the formal definition of the difference subgroup graph 
$D(G)$ of a group $G$. 

\begin{definition} Let $G$ be a finite group and $S$ be the collection of all non-trivial proper subgroups of $G$. The difference subgroup graph $D(G)$ of a group $G$ is defined to be a graph with $S$ as the set of vertices and two distinct vertices $H$ and $K$ are adjacent if $\langle H,K \rangle=G$ and $HK\neq G$. The deleted difference subgroup graph of $G$, denoted by $D^*(G)$, is defined as the graph obtained by removing the isolated vertices, if any, from $D(G)$.
\end{definition}

Note that while the definition above is stated for finite groups, it extends naturally to infinite groups without modification. Additionally, it follows directly from the definition that $D(G)$ is an undirected graph. Furthermore, it is straightforward to verify that for abelian groups—and more generally, for Dedekind groups (groups in which every subgroup is normal) — we have $\langle H,K\rangle=HK$. Consequently, the difference subgroup graph $D(G)$ is edgeless. 
We can extend this observation further: if $G$ is an Iwasawa group (a group in which every subgroup is permutable), then $D(G)$ is again edgeless. This raises a natural question: Do there exist non-Iwasawa groups $G$ for which $D(G)$ is still edgeless. The answer is yes. For example, $\mathbb{Z}_4 \times Q_8$ is not an Iwasawa group, yet its difference subgroup graph $D(G)$ is edgeless. We discuss more about this in Corollary \ref{edgeless-impliy-nilpotent}.

In this current paper, we first establish fundamental structural properties of $D(G)$ and $D^*(G)$ in Sections \ref{basic-section} and \ref{Dstar-section}. Section \ref{forbidden-section} presents some forbidden subgraph characterizations of $D(G)$. We then investigate how the independence number and clique number of $D(G)$ reflect solvability/nilpotency of the underlying group in Sections \ref{independence-section} and \ref{clique-section}. Finally, we conclude with a discussion of open problems and future studies.
 
\section{Basic Results}\label{basic-section}
We begin by recalling some well-known results from finite group theory, which will be used throughout the paper. 

\begin{proposition}\label{bunch}
The following holds for finite groups:
\begin{enumerate}
    \item (\cite{miller2}) A simple group of composite order cannot contain a maximal subgroup of prime order.
    \item (Ex. 7, Sec. 10.5, \cite{robinson-book}) A finite group with an abelian maximal subgroup is solvable with derived length at most $3$.
    \item (\cite{ore}) If $G$ is a solvable group and $M$ and $N$ are two maximal subgroups of $G$, then either $MN=G$ or $M$ and $N$ are conjugate in $G$.
    \item (Lemma 6, \cite{garonzi}) Let $G$ be a finite group such that $G=HK$ for two subgroups $H,K$ of $G$. Then $G=(xHx^{-1})(yKy^{-1})$ for all $x,y \in G$. 
\end{enumerate}
\end{proposition}

We now prove some basic properties related to $D(G)$.

\begin{proposition}\label{basic_results}
The following are true:
    \begin{enumerate}
        \item Let $G$ be a group and $H$ be a subgroup of $G$ such that $H\overline{H}=G$ for some conjugate $\overline{H}$ of $H$. Then $H=G$.
        \item If $M$ is a maximal subgroup of $G$ which is not normal in $G$, then for all $g\in G\setminus M$, we have $M\sim gMg^{-1}$ in $D(G)$.
        \item If $N$ is a non-trivial proper normal subgroup of $G$, then $N$ is an isolated vertex in $D(G)$.
        \item If $H \sim K$ in $D(G)$, then for all $g \in G$, $gHg^{-1}\sim gKg^{-1}$ in $D(G)$.
        \item Any vertex in $D(G)$ is either isolated or of degree $\geq 2$, i.e., there does not exist any leaf in $D(G)$.
        \item The degree of any non-isolated vertex in $D(G)$ is not unique.
        \item If $G\cong H \rtimes K$, then $D(K)$ is isomorphic to an induced subgraph of $D(G)$.
        \item If $N$ is a normal subgroup of $G$, then $D(G/N)$ is isomorphic to an induced subgraph of $D(G)$.
    \end{enumerate}
\end{proposition}
\pf \begin{enumerate}
    \item Let $g\in G$ such that $\overline{H}=gHg^{-1}$ and $H\overline{H}=G$. If $g \in H$, then clearly $H=G$. So we assume that $g \in G\setminus H$. Now, $$g^{-1}=(h_1)(gh_2g^{-1})\in H\overline{H}=G, \mbox{ i.e., } g=h^{-1}_1h^{-1}_2\in H, \mbox{ a contradiction.}$$
    \item As $M$ is maximal and not normal, we have $\langle M,gMg^{-1}\rangle=G$. Also by Proposition \ref{basic_results}(1), $M(gMg^{-1})\neq G$. Hence the lemma holds.
    \item It follows from the fact that $\langle N,H\rangle=NH$ for all subgroups $H$ of $G$.
    \item As $H \sim K$, we have $\langle H,K\rangle =G$ and $HK\neq G$. Let $x \in G$. Then $$x=h^{\alpha_1}_1k^{\beta_1}_1h^{\alpha_2}_2k^{\beta_2}_2\cdots h^{\alpha_r}_rk^{\beta_r}_r$$ for some $h_i\in H,k_i\in K$ and $\alpha_i,\beta_i \in \mathbb{Z}$. Thus $$gxg^{-1}=(gh^{\alpha_1}_1g^{-1})(gk^{\beta_1}_1g^{-1})\cdots (gh^{\alpha_r}_rg^{-1})(gk^{\beta_r}_rg^{-1})\in \langle gHg^{-1},gKg^{-1} \rangle=T\mbox{ (say)}.$$
Therefore $x \in g^{-1}Tg$ and hence $G \subseteq g^{-1}Tg$, i.e., $G\subseteq T$, i.e., $G=\langle gHg^{-1},gKg^{-1} \rangle$. Also, if $(gHg^{-1})(gKg^{-1})=G$, then $gHKg^{-1}=G$, i.e., $HK=G$, a contradiction. Thus $(gHg^{-1})(gKg^{-1})\neq G$. Hence $gHg^{-1}\sim gKg^{-1}$ in $D(G)$.
\item Let $H$ be a vertex in $D(G)$. If $H$ is isolated, then there is nothing to prove. If not, let $H\sim K$ in $D(G)$. Thus by Proposition \ref{bunch}(3), $H,K$ are not normal in $G$, i.e., $N_G(H),N_G(K)$ are proper subgroups of $G$. Choose $y \in N_G(H)\setminus N_G(K)$. Then by Proposition \ref{bunch}(4), we have $yHy^{-1}=H\sim yKy^{-1}\neq K$, i.e., $H$ has a neighbour other than $K$ in $D(G)$. Hence the result follows.
\item Let $H$ be a non-isolated vertex and $H\sim K$. Choose $y \in N_G(K)\setminus N_G(H)$. Then we have $H\sim K\sim yHy^{-1}\neq H$. As $deg(H)=deg(yHy^{-1})$, the result follows.
\item Let $K_1,K_2$ be two non-trivial proper subgroups of $K$ such that $K_1\sim K_2$ in $D(K)$. Then $\langle K_1,K_2\rangle=K$ and $K_1K_2\neq K$. Thus $\langle HK_1,HK_2\rangle=HK=G$ and $HK_1\cdot HK_2\neq G$. Note that $HK_1,HK_2$ are distinct, proper subgroups of $G$ and hence $HK_1\sim HK_2$ in $D(G)$. Similarly, it can be shown that $HK_1\sim HK_2$ in $D(G)$ implies that $K_1\sim K_2$ in $D(K)$.
\item Let $H/N,K/N$ be two non-trivial proper subgroups of $G/N$ such that $H/N\sim K/N$ in $D(G/N)$. Then $\langle H/N,K/N \rangle=G/N$ and $H/N \cdot K/N\neq G/N$, i.e., $\langle H,K\rangle/N=G/N$ and $HK/N\neq G/N$, i.e., $\langle H,K\rangle = G$ and $HK\neq G$, i.e., $H\sim K$ in $D(G)$. Similarly, it can be shown that if $H,K$ are two subgroups of $G$ containing $N$, then $H\sim K$ in $D(G)$ implies $H/N\sim K/N$ in $D(G/N)$.
\end{enumerate} \qed 

\begin{lemma}\label{min_no_of_edges}
    Let $H\sim K$ in $D(G)$. Then:
    \begin{enumerate}
        \item If $H$ and $K$ are conjugates, then $D(G)$ has at least $3$ edges.
        \item If $H$ and $K$ are not conjugates, then $D(G)$ has at least $4$ edges.
    \end{enumerate}
\end{lemma}
\pf It is clear that $H,K$ are not normal in $G$, i.e., $N_G(H),N_G(K)$ are proper subgroups of $G$. Also neither $N_G(H)$ nor $N_G(K)$ is contained in each other, as that would imply $H,K\subseteq N_G(H)$ or $N_G(K)$, i.e., $\langle H,K\rangle \neq G$. Moreover $G\neq N_G(H)\cup N_G(K)$. Let us choose $x \in G\setminus (N_G(H)\cup N_G(K))$, $y \in N_G(H)\setminus N_G(K)$ and $z \in N_G(K)\setminus N_G(H)$.

\begin{enumerate}
    \item As $H \sim K$, we have $H=yHy^{-1}\sim yKy^{-1}\neq K$ and $H\neq zHz^{-1}\sim zKz^{-1}=K$. Hence we get three distinct edges.
    \item Also $H\neq xHx^{-1}\sim xKx^{-1}\neq K$. As $H$ and $K$ are not conjugates, we have $H\neq xKx^{-1}$ and $K\neq xHx^{-1}$. Hence we get another edge distinct from the edge $H\sim K$.
\end{enumerate} \qed 

\begin{proposition}\label{basic_nilpotent}
    Let $G$ be a finite nilpotent group. Then:
    \begin{enumerate}
        \item No two conjugate subgroups are adjacent in $D(G)$.
        \item If $\langle H,K\rangle=G$, then $\langle xHx^{-1},yKy^{-1}\rangle=G$ for all $x,y \in G$.
        \item If $H\sim K$ in $D(G)$, then $xHx^{-1}\sim yKy^{-1}$ in $D(G)$ for all $x,y \in G$.
        \item If $D(G)$ has at least one edge, then $D(G)$ contains a $4$-cycle as an induced subgraph.
    \end{enumerate}
\end{proposition}
\pf \begin{enumerate}
    \item If possible, let $H$ and $gHg^{-1}$ be adjacent in $D(G)$. Clearly $H$ is not maximal in $G$, as $G$ is nilpotent, i.e., maximal subgroups are normal. So $H$ is properly contained in some normal maximal subgroup $M$ of $G$. As $M$ is normal, $gHg^{-1}\subseteq M$. So, we have $\langle H,gHg^{-1} \rangle\subseteq M\neq G$, i.e., $H\nsim gHg^{-1}$ in $D(G)$, a contradiction. 
    \item If $\langle xHx^{-1},yKy^{-1}\rangle\neq G$, then there exists a maximal subgroup $M$ of $G$ such that $\langle xHx^{-1},yKy^{-1}\rangle \leq M$, i.e.,      
    $xHx^{-1},yKy^{-1}\leq M$. As $M$ is normal in $G$, we have $H=x^{-1}(xHx^{-1})x\leq x^{-1}Mx=M$ and similarly $K\leq M$. Thus we have $\langle H,K\rangle \leq M\neq G$, a contradiction.
    \item As $H \sim K$ in $D(G)$, we have $\langle H,K\rangle=G$ and $HK\neq G$. So, by Proposition \ref{basic_nilpotent}(2), we have $\langle xHx^{-1},yKy^{-1}\rangle=G$. If $(xHx^{-1})(yKy^{-1})=G$, then by Proposition \ref{bunch}(4), we have $\{x^{-1}(xHx^{-1})x\}\{y^{-1}(yKy^{-1})y\}=G$, i.e., $HK=G$, a contradiction. Thus $(xHx^{-1})(yKy^{-1})\neq G$ and hence $xHx^{-1}\sim yKy^{-1}$ in $D(G)$.
    \item Let $H\sim K$ in $D(G)$. Then $H,K$ are non-normal subgroups of $G$ and by Proposition \ref{basic_nilpotent}(1), $H$ and $K$ are not conjugate. Let $x \in G\setminus (N_G(H)\cup N_G(K))$. Then $H,K,xHx^{-1},xKx^{-1}$ are four distinct subgroups of $G$. As $H\sim K$, by Proposition \ref{basic_nilpotent}(3), $H\sim xKx^{-1}$ and $K\sim xHx^{-1}$. Also, as $G$ is nilpotent, by Proposition \ref{basic_nilpotent}(1), we have $H\not\sim xHx^{-1}$ and $K\not\sim xKx^{-1}$. Thus $H\sim K\sim xHx^{-1}\sim xKx^{-1}\sim H$ is an induced $4$-cycle in $D(G)$.
\end{enumerate} \qed

\begin{theorem}\label{simple-connected}
    $D(G)$ is connected if and only if $G$ is simple.
\end{theorem}
\pf If $D(G)$ is connected, then it has no isolated vertices. Thus by Proposition \ref{basic_results}(3), $G$ has no proper normal subgroup, i.e., $G$ is simple.

Conversely, let $G$ be simple. Then, by Proposition \ref{basic_results}(2) any maximal subgroup is adjacent to all its conjugates in $D(G)$, i.e., a maximal subgroup and its conjugates form a complete subgraph in $D(G)$.

{\it Claim 1:} If $H$ is a non-trivial subgroup contained in some maximal subgroup $M$ of $G$, then $H$ is adjacent to some conjugate of $M$ in $D(G)$ and $d(H,M)=2$.\\
{\it Proof of Claim 1:} If $H$ is contained in all conjugates of $M$, then $H\subseteq \bigcap_{g\in G}gMg^{-1}=N$(say). As $H$ is non-trivial and $N$ is a normal subgroup of a simple group $G$, we get a contradiction. Thus there exists a conjugate $\overline{M}$ of $M$ such that $H\not\subseteq \overline{M}$. As $\overline{M}$ is a maximal subgroup, $\langle H,\overline{M}\rangle=G$ and $H\overline{M}\subseteq M\overline{M}\neq G$, we get $H\sim \overline{M}\sim M$ in $D(G)$. Hence Claim 1 holds.

Claim 1 shows that any non-trivial subgroup of $G$ is adjacent to some maximal subgroup of $G$ in $D(G)$. Let $H$ be a non-maximal subgroup of $G$. It suffices to show that there exists a path joining $H$ and any maximal subgroup of $G$ not containing $H$ in $D(G)$. Let $N$ be a maximal subgroup of $G$ not containing $H$. Let $M$ be a maximal subgroup containing $H$. If $MN\neq G$, we have $HN\neq G$ and  $\langle H,N\rangle=G$, i.e., $H \sim N$. So we assume that $MN=G$.

If $A=M\cap N$ is non-trivial, then by Claim 1, $A$ must be adjacent to some conjugate $M'$ of $M$ and some conjugate $N'$ of $N$. Thus we get a path joining $H$ and $N$ via $M',A$ and $N'$. So we assume that $M\cap N$  is trivial. Thus $H \cap N$ is also trivial.

So, we have $$|G|=|MN|=\dfrac{|M||N|}{|M\cap N|}>\dfrac{|H||N|}{|H \cap N|}=|HN|,$$
i.e., $HN\neq G$ and hence $H \sim N$. Thus in any case we get a path joining $H$ and $N$. \qed

\begin{theorem}\label{bipartite-implies-nilpotent}
    If $D(G)$ is triangle-free or bipartite, then $G$ is nilpotent.
\end{theorem}
\pf If possible, let $G$ be non-nilpotent. Then $G$ has a maximal subgroup $M$ which is not normal in $G$. Thus, the number of conjugates of $M$ in $G$ is $$[G:N_G(M)]=[G:M]\geq 3.$$ Hence $M$ has at least two other conjugates, say $M_1$ and $M_2$ in $G$. Now, by Proposition \ref{basic_results}(2), $M,M_1,M_2$ forms a triangle in $D(G)$, contradicting that $D(G)$ is bipartite.\qed 

\begin{corollary}\label{edgeless-impliy-nilpotent}
    If $D(G)$ is edgeless, then $G$ is nilpotent.
\end{corollary}
\pf If $G$ is non-nilpotent, then by above theorem, $G$ has a non-normal maximal subgroup and there exists a triangle in $D(G)$, contradicting that $D(G)$ is edgeless.\qed 

\begin{corollary}\label{girth}
    If $D(G)$ has at least one edge, then girth of $D(G)$ is $3$ or $4$.
\end{corollary}
\pf If $G$ is non-nilpotent, by Theorem \ref{bipartite-implies-nilpotent}, $G$ has a triangle. If $G$ is nilpotent, by Proposition \ref{basic_nilpotent}(4), $D(G)$ has an induced $4$-cycle. Hence the corollary follows.\qed 

\begin{remark}
    The converse of Theorem \ref{bipartite-implies-nilpotent} is not true. If $G\cong (\mathbb{Z}_2 \times \mathbb{Z}_2 \times \mathbb{Z}_2)\rtimes (\mathbb{Z}_2\times \mathbb{Z}_2)$ [GAP id: $(32,49)$], then $D(G)$ is not bipartite. Moreover, nilpotent groups may yield graphs with girth $3$, for example $G\cong (\mathbb{Z}_3 \times \mathbb{Z}_3 )\rtimes \mathbb{Z}_3$. 
\end{remark}

\begin{theorem}
    $D(G)$ cannot have a universal vertex, and hence $D(G)$ is never complete. 
\end{theorem}
\pf If $D(G)$ has a universal vertex, then $D(G)$ is connected and hence $G$ is simple. Let $H$ be a universal vertex in $D(G)$. Then $H$ is both a maximal and minimal subgroup of $G$, i.e., $H$ is a maximal subgroup of prime order in $G$. However this contradicts Proposition \ref{bunch}(1).\qed 

\begin{remark}
    It is to be noted that $D(G)$ can not be a cycle. It follows from Theorem \ref{simple-connected} and Theorem \ref{bipartite-implies-nilpotent}.
\end{remark}

\section{Properties of $D^*(G)$}\label{Dstar-section}
As established in the previous section, $D(G)$ frequently contains isolated vertices, particularly normal subgroups. To obtain a more informative graph, we now focus on $D^*(G)$, the graph obtained by removing all isolated vertices from $D(G)$. We start by noting that $D^*(G)$ cannot be a tree, as $D(G)$ has no leaf (Proposition \ref{bunch}(5)). Next, we investigate when $D^*(G)$ admits a universal vertex. 

\begin{theorem}\label{Dstar(G)_universal}
    If $D^*(G)$ has a universal vertex, then $G\cong \mathbb{Z}^\beta_q \rtimes \mathbb{Z}_{p^\alpha}$, where $p,q$ are distinct primes.
\end{theorem}
\pf Let $H$ be a universal vertex in $D^*(G)$. We break the proof down into several steps, each stated as a claim.

{\it Claim 1:} $G$ is not nilpotent.\\
{\it Proof of Claim 1:} Suppose $G$ is nilpotent and $H\sim K$ in $D(G)$. Hence by Proposition \ref{basic_results}(3), $H$ and $K$ are non-normal in $G$, i.e., $N_G(H)$ and $N_G(K)$ are proper subgroups of $G$. Thus $N_G(H)\cup N_G(K) \neq G$. We choose $g \in G \setminus (N_G(H)\cup N_G(K))$. Hence $gHg^{-1}\neq H$ and $gKg^{-1}\neq K$. Also, by Proposition \ref{basic_results}(4), we have $gHg^{-1}\sim gKg^{-1}$ in $D(G)$. Thus $gHg^{-1}$ is not isolated in $D^*(G)$ and hence it must be adjacent to $H$. As $G$ is nilpotent, any maximal subgroup of $G$ is normal in $G$. Thus $H$ is not a maximal subgroup of $G$, i.e., $H$ is properly contained in some maximal subgroup $M$ of $G$. Therefore $gHg^{-1}\subseteq gMg^{-1}=M$, i.e., both $H,gHg^{-1}\subseteq M$. However, this implies that $H \nsim gHg^{-1}$ in $D(G)$, a contradiction. Hence the claim holds.

{\it Claim 2:} $H$ is a maximal subgroup of $G$ and hence $N_G(H)=H$.\\
{\it Proof of Claim 2:} Suppose, $H$ is not maximal in $G$. Then there exist a maximal subgroup $M$ of $G$ containing $H$ properly. Therefore, $\langle H,M\rangle=M\neq G$ and hence $H\not\sim M$. Now, as $H$ is a universal vertex, $M$ must be an isolated vertex in $D(G)$. This implies $M\lhd G$ as otherwise $M$ is adjacent to its conjugates (by Proposition \ref{basic_results}(2)), thereby making $M$ to be a non-isolated vertex. Therefore $M$ contains $H$ and all its conjugates.

Choose $z \in N_G(K)\setminus N_G(H)$. Then, by Proposition \ref{basic_results}(4), $zKz^{-1}=K\sim zHz^{-1}\neq H$ in $D(G)$. Now, as $H$ is a universal vertex in $D^*(G)$, we have $H \sim zHz^{-1}$, i.e., $\langle H,zHz^{-1}\rangle=G$. However this contradicts that $H$ and its conjugates are contained in $M$. Hence Claim 2 holds.

{\it Claim 3:} $H$ is cyclic Sylow subgroup of $G$.\\
{\it Proof of Claim 3:} We first show that $H$ is cyclic $p$-group. If not, it should have at least two distinct maximal subgroups $H_1$ and $H_2$. Then $H=\langle H_1,H_2\rangle$. Let $\overline{H}$ be a conjugate of $H$ which is different from $H$. If both $H_1,H_2\subseteq \overline{H}$, then $\overline{H}=\langle H_1,H_2\rangle=H$, a contradiction. Thus, without loss of generality, let $H_1\not\subseteq \overline{H}$. Then $\langle \overline{H},H_1\rangle=G$ (as $\overline{H}$ is maximal in $G$). Also $\overline{H}H_1\subseteq \overline{H}H\neq G$. Thus $\overline{H}\sim H_1$ and $H_1$ is not isolated in $D(G)$. As $H$ is a universal vertex, we have $H\sim H_1$, a contradiction, as $H_1\subseteq H$. Thus $H$ has a unique maximal subgroup and hence $H$ is a cyclic $p$-group. As every $p$-group is contained in some Sylow $p$-subgroup and $H$ is a maximal subgroup of $G$, $H$ must be a cyclic Sylow $p$-subgroup of $G$. Thus Claim 3 holds.

From Claim 2 and 3 and Proposition \ref{bunch}(2), $G$ is solvable. 

{\it Claim 4:} $|G|$ has exactly two distinct prime factors, i.e., $|G|=p^\alpha q^\beta$.\\
{\it Proof of Claim 4:} As $H$ is a Sylow $p$-subgroup of $G$, let $|H|=p^\alpha$. As $G$ is solvable, every maximal subgroup is of prime power index, i.e., $[G:H]=q^\beta$, where $q$ is a prime. So, $|G|=p^\alpha q^\beta$. Again, as $G$ is not nilpotent, we have $q\neq p$. Thus Claim 4 holds.




{\it Claim 5:} $G$ has exactly one subgroup of order $p^i$ for $i=1,2,\ldots,\alpha-1$.\\
{\it Proof of Claim 5:} Suppose $G$ has more than one subgroup of order $p^i$ for some $i$, say $H_1$ and $H_2$. Since $H$ is a cyclic group of order $p^\alpha$, it can contain at most one of them. Suppose $H_1\not\subseteq H$. Also $H_1$ is contained in some conjugate $\overline{H}$ of $H$. Since $H$ is maximal, we have $\langle H_1,H\rangle=G$ and $HH_1\subseteq H\overline{H}\neq G$, i.e., $H \sim H_1$. Now, as $H$ is a universal vertex in $D^*(G)$ and any inner automorphism of $G$ induces an automorphism of $D(G)$,  $\overline{H}$ is also a universal vertex of $D^*(G)$. So, we must have $H_1\sim \overline{H}$. However this is a contradiction as $H_1\subseteq \overline{H}$. So Claim 5 holds.

As $G$ is solvable, $G$ must have a maximal normal subgroup $M$ of prime index. So $M$ must be of order $p^{\alpha-1}q^\beta$ or $p^\alpha q^{\beta-1}$. The latter cannot hold because $H$ is of order $p^\alpha$ and $H$ is itself a maximal subgroup of $G$. Thus $|M|=p^{\alpha-1}q^\beta$ and $M\lhd G$. 

{\it Claim 6:} Order of every element in $G\setminus M$ is a multiple of $p^\alpha$.\\
{\it Proof of Claim 6:} Firstly note that from Claim 5, it follows that all elements of order $p^i$ for $i=1,2,\ldots,\alpha-1$ belong to $M$. Let $Q$ be a Sylow $q$-subgroup of $M$. Note that $Q$ is also a Sylow $q$-subgroup of $G$. As Sylow $q$-subgroups of $G$ are conjugate in $G$ and $M$ is normal in $G$, all Sylow $q$-subgroups of $G$, and thereby all $q$-subgroups of $G$ are contained in $M$. Let $x \in G\setminus M$. If possible, let $\circ(x)=p^i q^j$ where $i<\alpha-1$. Set $X=\langle x \rangle$. Then $MX=G$, as $M$ is maximal and normal in $G$. Note that the unique subgroup of order $p^i$ (as shown in Claim 5) is contained in both $M$ and $X$. So, we get $$p^\alpha q^\beta=|G|=|MX|=\dfrac{|M||X|}{|M\cap X|}=\dfrac{p^{\alpha-1} q^\beta \cdot p^iq^j}{p^i\cdot q^t}=p^{\alpha-1}q^{\beta+j-t}, \mbox{ a contradiction.}$$ Thus Claim 6 holds.

{\it Claim 7:} Sylow $q$-subgroup is normal in $G$.\\
{\it Proof of Claim 7:} Recall that $Q$ is a Sylow $q$-subgroup of $M$ and $M\lhd G$. Then by Frattini's argument, $G=N_G(Q)M$. Thus $N_G(Q)\not\subseteq M$ and $q^\beta \mid |N_G(Q)|$. By Claim 6, $N_G(Q)$ contains an element of order a multiple of $p^\alpha$, i.e., $p^\alpha\mid  |N_G(Q)|$. Thus $p^\alpha q^\beta \mid |N_G(Q)|$, i.e., $N_G(Q)=G$, i.e., $Q \lhd G$.

{\it Claim 8:} $Q$ is elementary abelian $q$-group.\\
{\it Proof of Claim 8:} As $Q\lhd G$, $Q$ contains a minimal normal subgroup of $G$. Again, as $G$ is solvable, a minimal normal subgroup is isomorphic to $\mathbb{Z}^t_q$. If $t<\beta$, then we get a subgroup of order $p^\alpha q^t$ containing $H$ in $G$, a contradiction to maximality of $H$ in $G$. Thus $t=\beta$, i.e., $Q\cong \mathbb{Z}^\beta_q$.

Combining all the claims, we get $G\cong \mathbb{Z}^\beta_q \rtimes \mathbb{Z}_{p^\alpha}$. \qed

\begin{corollary}\label{D^*(G)-complete}
    If $D^*(G)$ is complete, then $G\cong \mathbb{Z}_q \rtimes \mathbb{Z}_{p^\alpha}$.
\end{corollary}
\pf As $D^*(G)$ is complete, all of its vertices are universal and hence from the above theorem, we get $G\cong \mathbb{Z}^\beta_q \rtimes \mathbb{Z}_{p^\alpha}$. We recall that Sylow $p$-subgroup $H$ is a maximal subgroup. Now, if $\beta>1$, then $G$ has a subgroup $K$ of order $q$. Then $\langle H,K\rangle=G$ and $|HK|=qp^\alpha<|G|$, i.e., $HK\neq G$. So, we have $H\sim K$. As $D^*(G)$ is complete, $K$ must also be a universal vertex in $D^*(G)$. But this implies that $K$ must be a maximal subgroup of $G$ (using Claim 2 of previous theorem). However, this cannot be true as $K$ is properly contained in some Sylow $q$-subgroup of $G$ (as $\beta>1$). So, we have $\beta=1$ and the corollary follows. \qed 

\begin{remark}\label{complete-no-of-vertices}
    If $G\cong \mathbb{Z}_q \rtimes \mathbb{Z}_{p^\alpha}$, then $D^*(G)$ is complete and the number of vertices in $D^*(G)$ is $n_p$, the number of Sylow $p$-subgroups of $G$. So we have $1\neq n_p=1+pl=q$.
\end{remark}

\begin{theorem}\label{not-a-cycle}
    If $D^*(G)$ is a cycle, it has length $3$ or $4$.
\end{theorem}
\pf We start by noting that $D^*(S_3)\cong C_3$ and $D^*(D_4)\cong C_4$. Now, if possible, let $D^*(G)$ be a $l$-cycle with $l\geq 5$. As $D^*(G)$ is triangle-free, by Theorem \ref{bipartite-implies-nilpotent}, $G$ is nilpotent. Hence by Proposition \ref{basic_nilpotent}(4), $D^*(G)$ must have an induced $4$-cycle, a contradiction.\qed

\section{Forbidden Subgraph Characterization of $D(G)$}\label{forbidden-section}
In Theorem \ref{bipartite-implies-nilpotent}, we established that bipartiteness (the absence of odd cycles) in $D(G)$ forces nilpotency of the underlying group. In this section, we explore further forbidden subgraph characterizations, specifically examining when $D(G)$ is clawfree or cograph, and the implications for the structure of $G$.
\begin{theorem}\label{clawfree_implies_supersolvable}
    Let $G$ be a group such that $D(G)$ is clawfree. Then $G$ is supersolvable.
\end{theorem}
\pf Suppose it is not true and $G$ be the minimum example, i.e., $G$ be the group of minimum order which is not supersolvable but $D(G)$ is clawfree. Since $G$ is not supersolvable, there exists a maximal subgroup $M$ of $G$ such that $[G:M]$ is not prime. Clearly this implies that $M$ is not normal in $G$. Let $M_1$ be a distinct conjugate of $M$ in $G$. 

If $M$ is not cyclic, choose $x \in M\setminus M_1$ and set $H=\langle x \rangle$. Then $H$ is a proper non-trivial subgroup of $G$ (as $M$ is not cyclic). Again as $H$ and $m\cap M_1$ are proper subgroups of $M$, we have $M\neq H \cup (M \cap M_1)$. Choose $y \in M\setminus [H \cup (M \cap M_1)]$ and set $K=\langle y \rangle$. Then $M,M_1,H,K$ form a claw with $M_1$ being the vertex of degree $3$, a contradiction. 

Thus $M$ is cyclic. Also as $G$ has a maximal subgroup $M$ which is cyclic, $G$ is solvable. Thus the index $[G:M]=q^k$ where $q$ is a prime and $k>1$. Note that, except $q$, all other Sylow subgroups are cyclic. Moreover Sylow $q$-subgroups are not cyclic, as otherwise $G$ will be supersolvable. Let $Q$ be a non-cyclic Sylow $q$-subgroup of $G$.

If $q\neq 2$, i.e., $q$ is an odd prime $\geq 3$, $Q$ has at least $1+q\geq 4$ subgroups of order $q$. As $M$ is cyclic, among these subgroups of order $q$, at most one can be contained in $M$, i.e., there exist at least three subgroups $Q_1,Q_2,Q_3$ of $Q$ of order $q$ which are not contained in $M$. Clearly $\langle M,Q_i\rangle=G$ for $i=1,2,3$ and $\langle Q_i,Q_j\rangle \subseteq Q\neq G$. If $MQ_i=G$, then we have $$q^2\leq q^k=\dfrac{|G|}{|M|}=\dfrac{|Q_i|}{|M\cap Q_i|}\leq q,\mbox{ a contradiction.}$$ Thus $MQ_i\neq G$ and hence $M,Q_1,Q_2,Q_3$ form a claw in $D(G)$ with $M$ being the vertex of degree $3$, a contradiction. Hence $q=2$ and $[G:M]=2^k$.

{\it Case 1:} $|M|$ is odd. As $Q$ is a non-cyclic Sylow $2$-subgroup of $G$, $Q$ has at least three proper non-trivial subgroups $Q_1,Q_2,Q_3$. As $|M|$ is odd, we have $Q_i \not\subseteq M$. Hence $M,Q_1,Q_2,Q_3$ form a claw in $D(G)$ with $M$ being the vertex of degree $3$, a contradiction.

{\it Case 2:} $|M|$ is even. As $M$ is cyclic, there exists a unique element $z \in M$ such that $\circ(z)=2$. Also, all conjugates of $M$ contains a unique element of order $2$. Since $[G:M]=2^k\geq 4$, $M$ has at least three distinct conjugates $M_1,M_2,M_3$ other than $M$. Let $M_i=g_i M g^{-1}_i$ and $z_i=g_i z g^{-1}_i\in M_i$ for $i=1,2,3$. 

If any one of these $z_i$'s belongs to $M$, say $z_1\in M$, then $z_1=g_1zg^{-1}_1=z$, i.e., $g_1z=zg_1$, i.e., $g_1\in C(z)$, the centralizer of $z$ in $G$. Also as $M$ is abelian and $z \in M$, we have $M\subseteq C(z)$. Thus, as $g_1\not\in M$, we have $G=\langle M,g_1\rangle=C(z)$, i.e., $z \in Z(G)$, the center of $G$. Let $N=\langle z \rangle$. Then $N$ is normal in $G$.

Consider the group $G/N$. Then $M/N$ is maximal subgroup of $G/N$ which is cyclic and non-normal in $G/N$. Also $[G/N:M/N]=2^k$. Thus $G/N$ is a non-supersolvable group of order $|G|/|N|$. Thus minimality of $G$, $D(G/N)$ has a claw, say $A/N,B_1/N,B_2/N,B_3/N$ where $A/N$ is the only degree of vertex $3$. Thus $A,B_1,B_2,B_3$ form a claw in $D(G)$, a contradiction. So, we must have $z_1,z_2,z_3 \not\in M$. Clearly $M\sim \langle z_i \rangle$ in $D(G)$. As $\langle z_i,z_j\rangle\neq G$ (which we prove in the next claim), we get a claw $M,\langle z_1\rangle,\langle z_2\rangle,\langle z_3\rangle$, a contradiction. Thus to be prove the theorem, it suffices to prove the following claim.

{\it Claim:} $\langle z_i,z_j\rangle\neq G$. As $G$ is solvable, $G$ has maximal subgroup $L$ such that $L\lhd G$ and $[G:L]=p$, a prime. If $p\neq 2$, then $L$ contains all Sylow $2$-subgroups of $G$. In particular $z_i,z_j \in L$, i.e., $\langle z_i,z_j\rangle\subseteq L \neq G$. If $p=2$, we have $[G:M]=2^k$ and $[G:L]=2$ with $k>1$. Note that at least one Sylow $p$-subgroup of $G$ is not normal for $p\neq 2$, because if all odd Sylow subgroups, say $P_1,P_2,\ldots,P_t$ of $G$ are normal, and as all of them are cyclic, then $P=P_1P_2\cdots P_t$ is a cyclic normal subgroup of $G$ and $G/P$ is a $2$-group, which implies $G$ is supersolvable, a contradiction. Without loss of generality, let $P_1$ be a cyclic odd Sylow subgroup of $G$ which is not normal in $G$. Thus $P_1$ has at least four conjugates, including $P_1$, in $G$ and as $L$ is normal in $G$, all of these conjugates are contained in $L$. Also exactly one of these four is contained in $M$ and the other three form a claw with $M$, a contradiction. Thus $p=2$ cannot hold. This proves the claim and the theorem. \qed

It is to be noted that clawfreeness of $D(G)$ does not imply nilpotency of $G$, as $D(S_3)$ is clawfree.

\begin{theorem}\label{cograph_implies_solvable}
    Let $G$ be a group such that $D(G)$ is cograph. Then $G$ is solvable.
\end{theorem}
\pf It is enough to prove the following: Given a finite non-solvable group $G$, there exist two conjugate maximal subgroups $M_1$ and $M_2$ of $G$ such that there exist two elements $a\in M_1\setminus M_2$, and $b \in M_2\setminus M_1$ such that the subgroup $\langle a,b\rangle$ generated by $a$ and $b$ is a proper subgroup of $G$. (because in this case we get $\langle b \rangle\sim M_1\sim M_2\sim \langle a\rangle$ as an induced $P_4$ in $D(G)$)

We prove the above statement by induction. Suppose it holds for all non-solvable groups of order $<n$ and let $G$ be a non-solvable group of order $n$.

{\it Claim 1:} For any non-trivial normal subgroup $K$ of $G$, $G/K$ is solvable.\\
{\it Proof of Claim 1:} If $G/K$ is non-solvable, by induction hypothesis, $G/K$ has two maximal subgroups (conjugates) $M_1/K$ and $M_2/K$ such that $aK \in M_1/K \setminus M_2/K$ and $bK \in M_2/K \setminus M_1/K$ such that $\langle aK,bK\rangle \neq G/K$. Then $M_1$ and $M_2$ are the conjugate, maximal subgroups of $G$. Note that since $K\leq M_1\cap M_2$, we have $a \in M_1\setminus M_2$ and $b \in M_2 \setminus M_1$. Thus Claim 1 holds.

{\it Claim 2:} $G$ has a unique minimal normal subgroup.\\
{\it Proof of Claim 2:} If $K$ and $L$ are two distinct minimal normal subgroups of $G$, then $K\cap L$ is trivial, and $G\cong G/(K\cap L)$ is embedded as a subgroup of $G/K\times G/L$. As both the quotients appearing in the product are solvable (by Claim 1), $G$ is solvable, a contradiction. Thus Claim 2 holds.

Let $K$ be the unique minimal normal subgroup of $G$. Moreover, $K$ is not solvable and hence of even order. Let $T$ be a Sylow $2$-subgroup of $K$. Then by Frattini argument, $G=K\cdot N_G(T)$.

{\it Claim 3:} $N_G(T)$ is a proper subgroup of $G$.\\
{\it Proof of Claim 3:} If $N_G(T)=G$, then $T$ is normal in $G$, i.e., $T\lhd K$. Now $K/T$ being an odd-order group is solvable and $T$ being a $2$-group is solvable, thereby making $K$ solvable, a contradiction. Thus Claim 3 holds.

Let $M_1$ be a maximal subgroup of $G$ containing $N_G(T)$ and $S$ be a Sylow $2$-subgroup of $G$ containing $T$. Then $T=S\cap K \lhd N_G(S)$. Therefore $N_G(S)\leq N_G(T)\leq M_1$. Hence by using Frattini argument, one can prove that $M_1$ is self-normalizing, i.e., $M_1$ is not normal in $G$. 

Note that $K$ is not contained in $M_1$, as otherwise $G=KN_G(T)\leq M_1$. Since $M_1$ does not contain the unique minimal normal subgroup $K$ of $G$, we see that $M_1$ contains no non-trivial normal subgroup of $G$ at all.

Let $M_2$ be a conjugate of $M_1$ in $G$ and $\mathcal{I}(M_1)$ be the set of all involutions of $M_1$. Note that $\mathcal{I}(M_1)$ is non-empty as $T\leq M_1$.

{\it Claim 4:} $\mathcal{I}(M_1)$ is not contained in $\mathcal{I}(M_2)$.\\
{\it Proof of Claim 4:} If it were, we would have $\mathcal{I}(M_1)=\mathcal{I}(M_2)$ since $\mathcal{I}(M_1)$ and $\mathcal{I}(M_2)$ have equal cardinality (recall that $M_1$ and $M_2$ are conjugates in $G$). But in that case, we have $\langle \mathcal{I}(M_1)\rangle=\langle \mathcal{I}(M_2) \rangle \lhd \langle M_1,M_2\rangle = G$, contrary to the fact that $M_1$ contains no non-trivial normal subgroup of $G$. Thus Claim 4 holds.

Similarly, we can show that $\mathcal{I}(M_2)$ is also not contained in $\mathcal{I}(M_1)$. So, we may choose an involution $a\in M_1\setminus M_2$ and an involution $b\in M_2\setminus M_1$. And the subgroup generated by them is dihedral, i.e., solvable. Hence $\langle a,b\rangle\neq G$.\qed

Note that $D(G)$ may be cograph with $G$ being non-supersolvable, e.g., $A_4$.

\section{Independence number of $D(G)$}\label{independence-section}
In this section, we investigate how the independence number $\alpha(D(G))$ — the size of the largest set of pairwise non-adjacent vertices—relates to the structural properties of $G$. We establish bounds that force $G$ to be non-nilpotent, supersolvable, or solvable.
\begin{theorem}
    Let $G$ be a finite group such that $D(G)$ has at least one edge and independence number of $D(G)\leq 5$. Then $G$ is not nilpotent.
\end{theorem}
\pf Suppose on the contrary, let $G$ be nilpotent. Then all maximal subgroups of $G$ are normal in $G$ and hence are isolated vertices in $D(G)$ and hence contribute to the independence number of $D(G)$. Thus $G$ has at most $5$ maximal subgroups. Also, if number of maximal subgroups of $G$ is $\leq 2$, then $G$ is cyclic and hence $D(G)$ is edgeless, a contradiction. 

If $G$ has exactly three maximal subgroups and $G$ is non-cyclic, $G$ must be a $2$-group and for any maximal subgroup $M$ of $G$, we must have $Sub(M)\leq 4$ and hence $M$ is cyclic. Thus all proper subgroups of $G$ are cyclic and $G$ is a minimal non-cyclic group. Therefore by Miller-Moreno result \cite{miller-moreno} and using the fact that $G$ is a $2$-group, $G$ is isomorphic to $\mathbb{Z}_2\times \mathbb{Z}_2$ or $Q_8$. However, in both these cases, $D(G)$ is edgeless, a contradiction.  

If $G$ has exactly four maximal subgroups, then for each maximal subgroup $M$ of $G$, we must have $Sub(M)\leq 3$, i.e., $M$ must be cyclic. Arguing similarly, $G$ is a minimal non-cyclic nilpotent group with exactly $4$ maximal subgroups and by Miller-Moreno result, $G\cong \mathbb{Z}_3\times \mathbb{Z}_3$. But as $G$ is abelian, we must have $D(G)$ to be edgeless, a contradiction.

If $G$ has exactly five maximal subgroups, then each maximal subgroup is of prime order and by Miller-Moreno result, no such group exists, a contradiction. This proves the theorem. \qed

\begin{remark}
    The above bound is tight as independence number of $D(D_4)$ is $6$.
\end{remark}

\begin{theorem}
    Let $G$ be a finite group such that $D(G)$ has at least one edge and independence number of $D(G)\leq 13$. Then $G$ is either a $p$-group or $G$ is not nilpotent.
\end{theorem}
\pf Let $G$ be a group such that $D(G)$ has at least one edge and $\alpha(D(G))\leq 13$. If $G$ is a $p$-group, the theorem holds. Suppose $G$ is not a $p$-group. As $D(G)$ is not edgeless, $G$ is not cyclic. Again as $G$ is not a $p$-group, $G$ has at least four distinct maximal subgroups, $M_1,M_2,M_3,M_4$, say.We show that in this case, $G$ must be non-nilpotent. On the contrary, let $G$ be a nilpotent group which is not a $p$-group. Let $p_1,p_2,\ldots,p_k$ be the distinct prime factors of $|G|$ with $k\geq 2$.

{\it Claim:} $k=2$.\\
{\it Proof of Claim:} Suppose $k\geq 3$. Then $G$ has $k$ Sylow subgroups $P_1,P_2,\ldots,P_k$, say. Note that as $G$ is nilpotent and $|G|$ has at least $3$ distinct prime factors, the maximal subgroups and Sylow subgroups of $G$ are normal in $G$, and no Sylow subgroup is maximal in $G$. Thus we get at least $7$ isolated vertices in $D(G)$. Thus we must have $Sub(M_i)\leq 8$ for all $i$ (because otherwise proper non-trivial subgroups of $M_i$ along with the $7$ isolated vertices form an independent set of size $\geq 14$ in $D(G)$, a contradiction.) Now, as $M_i$'s are nilpotent, from Table 1 in \cite{das-mandal-number-of-subgroups}, each $M_i$ is either a $p$-group or cyclic. As $[G:M_i]$ is prime and $G$ has at least $3$ prime factors, $M_i$ can not be a $p$-group. So $M_i$'s are cyclic for all $i$, i.e., $G$ is a non-cyclic group such that all its maximal subgroups are cyclic, i.e., $G$ is a minimal non-cyclic group. However, by the classification of Miller-Moreno, this implies that $|G|$ has at most two prime factors, a contradiction. Hence Claim holds.

From the above claim, it follows that $G\cong P_1\times P_2$ where $P_i$'s are unique Sylow $p_i$-subgroups of $G$. Again, as $D(G)$ is not edgeless, $G$ is non-abelian. As a result, at least one of the two Sylow subgroups, say $P_1$, is non-abelian. Thus $|P_1|\geq p^3_1$ and $$Sub(P_1) \mbox{ is } \left\lbrace \begin{array}{ll}
=6, & \mbox{ if }P_1\cong Q_8 \\
\geq 10, & \mbox{ otherwise.} 
\end{array} \right.$$

If $Sub(P_1)\geq 10$, then we get an independent set of size $$9+3+1+1=14, \mbox{ a contradiction.}$$ where $9$ denotes the non-trivial subgroups of $P_1$ (including $P_1$), $3$ denotes the maximal subgroups of $G$ (other than possible $P_1$), $1$ stands for $P_2$ and the last $1$ denotes the subgroup $Z(P_1)\times P_2$ (note that this is not a maximal subgroup).

So, we must have $G\cong Q_8\times P_2$. If $P_2$ is abelian, then all subgroups of $G$ are normal in $G$ and hence all vertices in $D(G)$ are isolated, which is a contradiction. Thus $P_2$ must be non-abelian and we must have $Sub(P_2)\geq 10$. Again, we get a contradiction proceeding as in case of $Sub(P_1)\geq 10$. Hence the theorem holds.\qed

\begin{remark}
    The above bound is tight as independence number of $D(\mathbb{Z}_p\times D_4)$ is $14$, where $p$ is an odd prime.
\end{remark}

\begin{theorem}
    Let $G$ be a finite group such that $D(G)$ has at least one edge and independence number of $D(G)\leq 3$. Then $G$ is supersolvable.
\end{theorem}
\pf From the above theorem, it follows that $G$ is either a $p$-group or non-nilpotent. If it is $p$-group, it is supersolvable. Suppose $G$ is not a $p$-group and let $M$ be a maximal subgroup of $G$. Then $Sub(M)\leq 4$, as otherwise $M$ and its non-trivial subgroups form an independent set of size $\geq 4$ in $D(G)$. However, this implies that $M$ is cyclic, i.e., every maximal subgroup of $G$ is cyclic which implies that every Sylow subgroup of $G$ is cyclic, which in turn implies $G$ is supersolvable.  \qed

\begin{remark}
    The above bound is tight as independence number of $D(A_4)$ is $4$.
\end{remark}

\begin{theorem}
    Let $G$ be a finite group such that $D(G)$ has at least one edge and independence number of $D(G)\leq 14$. Then $G$ is solvable.
\end{theorem}
\pf Suppose $G$ be a minimal counterexample to this claim, i.e., $G$ is a smallest order non-solvable group such that $D(G)$ is not edgeless and $\alpha(D(G))\leq 14$. 

Let $N$ be any non-trivial normal subgroup of $G$. If $G/N$ is non-solvable, then $D(G/N)$ is not edgeless and $|G/N|<|G|$. So, by minimality of $G$, we must have $\alpha(D(G/N))\geq 15$. As $D(G/N)$ is an induced subgraph of $D(G)$, we must have $\alpha(D(G))\geq 15$, a contradiction. So $G/N$ is solvable and hence $N$ is not solvable. So, we must have $Sub(N)\geq 59$ (by Theorem 2.1 \cite{das-mandal-number-of-subgroups}). However, this implies that non-trivial subgroups of $N$ form an independent set of size at least $58$ in $D(G)$, a contradiction. So $G$ must be simple.

The above argument also shows that any proper subgroup of $G$ is solvable, i.e., $G$ is a finite minimal simple group. Thus by classification of finite minimal simple groups, $G$ is one of the following:
\begin{itemize}
    \item $PSL(2,2^p)$, where $p$ is a prime.
    \item $PSL(2,3^p)$, where $p$ is an odd prime.
    \item $PSL(2,p)$, where $p>3$ is a prime such that $5|p^2+1$.
    \item $PSL(3,3)$.
    \item The Suzuki group $Sz(2^{p})={}^{2}B_{2}(2^{p})$, where $p$ is an odd prime.
\end{itemize}
We show that the last one can not hold and for the rest of them $\alpha(D(G))\geq 15$. 

If $q=2^r$, then a Sylow $2$-subgroup $S$ of $PSL(2,q)$ is isomorphic to $\mathbb{Z}^r_2$. So $Sub(S)\geq 16$ for $r\geq3$. Thus non-trivial subgroups of $S$ form an independent set of size $\geq 15$ in $D(G)$ where $G\cong PSL(2,2^p)$, except when $p=2$. However one can compute and show that $\alpha(D(PSL(2,4))=15$.

If $q=p^r$, where $p$ is an odd prime, then $PSL(2,q)$ has a subgroup $H$ isomorphic to a dihedral group $D_{(q+1)/2}$ of order $q+1$. Unless $q=7$ or $13$, one can check that $Sub(H)\geq 16$, a contradiction. One can also separately check that $\alpha(D(PSL(2,7))=29$ and $\alpha(D(PSL(2,13)))=91$. 

As $PSL(3,3)$ has a subgroup isomorphic to $S_4$ and $Sub(S_4)=30$, we must have $\alpha(D(PSL(3,3)))\geq 29$. 

Coming to the last case, i.e., the Suzuki group, it is known that these are the only simple groups whose order is not divisible by $3$. However, we show that, in our case, $3$ must divide $|G|$.

As $G$ is non-solvable, at least one maximal subgroup of $G$, say $M$, must be non-supersolvable. Thus $Sub(M)=10$ or $15$ or $\geq 20$ (See Section 1.1 in \cite{das-mandal-number-of-subgroups}). If $Sub(M)\geq 20$, all non-trivial subgroups of $M$ forms an independent subset of size $\geq 19$ in $D(G)$, a contradiction. So $Sub(M)=10$ or $15$, and that implies $M\cong A_4$ or $SL(2,3)$ or $(\mathbb{Z}_2\times \mathbb{Z}_2)\rtimes \mathbb{Z}_9$. However, in each case, $3$ divides $|M|$ and hence $|G|$.\qed

\begin{remark}
    The above bound is tight as independence number of $D(A_5)$ is $15$.
\end{remark}

\section{Cliques in $D(G)$}\label{clique-section}
In this section, we investigate how the clique number $\omega(D(G))$ — the size of the largest complete subgraph — relates to the structural properties of $G$. We show that small clique numbers enforce supersolvability and solvability of the underlying group.
\begin{theorem}\label{omega<5_implies_supersolvable}
    Let $G$ be a group such that $\omega(D(G))\leq 4$. Then $G$ is supersolvable.
\end{theorem}
\pf Suppose there exists a group $G$ which is not supersolvable with $\omega(D(G))\leq 4$. As $G$ is not supersolvable, there exists a maximal subgroup $M$ of $G$ such that $[G:M]$ is not a prime. Thus $[G:M]=4$ or $[G:M]\geq 6$. Moreover as $[G:M]$ is not prime, $M$ is not normal in $G$. If $[G:M]\geq 6$, then $M$ and its other $5$ conjugates form a clique of size $\geq 6$, a contradiction. So, we must have $[G:M]=4$. 

Let the four conjugates of $M$ be $M=M_1,M_2,M_3$ and $M_4$. Note that $G \neq M_1\cup M_2\cup M_3\cup M_4$. If there exist any element $x \in G\setminus (M_1\cup M_2\cup M_3\cup M_4)$ such that $\circ(x)$ is odd, we set $H=\langle x\rangle$. Note that $\langle M_i,H\rangle=G$ and $$|HM_i|=\dfrac{|H||M_i|}{|H\cap M_i|}=\dfrac{|G|\cdot \circ(x)}{4\cdot |H\cap M_i|}\neq |G|, \mbox{ as }\circ(x) \mbox{ is odd.}$$ Thus $H$ along with the four conjugate maximal subgroups form a clique of size $5$, a contradiction. Thus every element of $G\setminus (M_1\cup M_2\cup M_3\cup M_4)$ must be of even order. 

Let $x$ be an element of minimum order in $G\setminus (M_1\cup M_2\cup M_3\cup M_4)$ and $H=\langle x \rangle$. Clearly $\langle M_i,H\rangle=G$ for all $i$. Again, as $\circ(x)$ is even, $\circ(x^2)<\circ(x)$ and hence $x^2 \in M_j$ for some $j \in \{1,2,3,4\}$. Thus for that $j$, $|H|/|H\cap M_j|=2$. Hence $$|HM_j|=|M_j| \cdot \dfrac{|H|}{|H\cap M_j|}=\dfrac{|G|}{4}\cdot 2=\dfrac{|G|}{2}\neq |G|,$$ i.e., $HM_j \neq G$. Now, by Proposition \ref{bunch} (4), it follows that $HM_i\neq G$ for all $i$. Thus $H$ along with four conjugate maximal subgroups form a clique of size $5$, a contradiction. Hence the theorem holds.

\begin{remark}\label{clique-remark}
    The bound in the above theorem is tight, as we have $\omega(D(A_4))=5$. Moreover, a similar result also holds for solvability, i.e., $\omega(D(G))\leq 7$ implies that $G$ is solvable. The proof follows on the same line of argument as above and uses the fact that a finite non-solvable group $G$ admits a maximal subgroup $M$ such that $[G:M]$ is neither a prime nor a prime-squared. 
\end{remark}

\section{Conclusion and open issues}
In this paper, we have introduced and systematically studied the difference subgroup graph $D(G)$ and its reduced version $D^*(G)$. Our investigation has revealed deep connections between graph-theoretic properties of $D(G)$ and fundamental group-theoretic properties of $G$. Despite the progress made, several intriguing questions remain open:
\begin{enumerate}
    \item \textbf{Connectivity of $D^*(G)$}: Although we have characterized when $D(G)$ is connected (Theorem \ref{simple-connected}), it remains an open issue when $D^*(G)$ is connected. We believe that the following is true and leave it as an open issue: {\it If $G$ is non-nilpotent, then $D^*(G)$ is connected.} Note that non-nilpotency is necessary, as if $G=GAP(32,49)$, then $D^*(G)$ is disconnected.
    \item \textbf{Girth in Nilpotent Groups}: It was proved in Theorem \ref{bipartite-implies-nilpotent}, that if $G$ is non-nilpotent, then the girth of $D(G)$ is $3$. Although we do not have any analogous result for nilpotent groups, we strongly suspect that the following is true: {\it If $G$ is a $p$-group, where $p$ is an odd prime and $D(G)$ has at least one edge, then $D(G)$ has girth $3$.}
    \item \textbf{Graph Isomorphism and Group Structure}: In light of Theorem \ref{simple-connected} and the fact that $A_5$ is the only simple group with $57$ non-trivial proper subgroups, it follows that $A_5$ is uniquely identifiable from its graph, i.e., $D(G)\cong D(A_5)$ implies $G\cong A_5$. However, it is not always true. For example, $D(S_3\times \mathbb{Z}_p) \cong D(S_3\times \mathbb{Z}_q)$, where $p,q$ are distinct odd primes (both having $9$ edges). Thus non-isomorphic groups can have isomorphic difference graphs. Similar examples are $D_4\times \mathbb{Z}_p$ and $Q_8\times \mathbb{Z}_p$. Another such pair is $D_5\times \mathbb{Z}_3$ and $\mathbb{Z}_5 \rtimes \mathbb{Z}_8$ both having $30$ edges. These, along with computational evidence, make us believe the following:
    \begin{enumerate}
        \item {\it Let $G$ and $H$ be finite groups such that $D(G)\cong D(H)$ and they are connected. Then $G\cong H$.}
        \item {\it Let $G$ and $H$ be finite groups such that $D(G)\cong D(H)$ and $G$ is nilpotent. Then $H$ is nilpotent.}
    \end{enumerate}

    \item \textbf{Additional Forbidden Subgraph Characterizations}: In Section \ref{forbidden-section}, we characterized various forbidden subgraphs in $D(G)$. One more in that list which we leave as an open issue is: {\it If $D(G)$ is perfect, then $G$ is solvable.}
    \item \textbf{Optimal Clique Number Bound}: In Remark \ref{clique-remark}, we mentioned that if the clique number is less than $8$, the underlying group is solvable. However, we think a better bound is possible, which we leave as a topic of further research: {\it $\omega(D(G))\leq 15$ implies $G$ is solvable.} 
\end{enumerate}

\section*{Acknowledgement}
The authors are grateful to Professor Geoffrey R. Robinson for providing fruitful insight towards the proof of Theorem \ref{cograph_implies_solvable}. The authors acknowledge the departmental funding of DST-FIST Sanction no. $SR/FST/MS-I/2019/41$. The first author is additionally supported by DST-SERB-MATRICS Sanction no. $MTR/2022/000020$, Govt. of India. The third author is supported by UGC PhD Junior Research Fellowship, Govt. of India.

\subsection*{Data Availability Statements}
Data sharing not applicable to this article as no datasets were generated or analysed during the current study.

\subsection*{Competing Interests} The authors have no competing interests to declare that are relevant to the content of this article.

\end{document}